\documentclass{article}
\usepackage{amsmath, amsfonts}
\usepackage{eufrak}
\usepackage{amssymb}
\usepackage{amscd}
\usepackage{exscale}

\begin{document}

\title{On an Algebraical computation of the tensor and the curvature for 3-Webs}

\author{Thomas B. Bouetou\footnote{currently at UMR 5030 (CNRS), D\'epartement des Sciences Math\'ematiques Universit\'e Montpellier II Case courrier 051-Place Eug\'ene Bataillon 34095 Montpellier CEDEX 05, France e-mail:tbouetou@darboux.math.univ-montp2.fr}\\ \'Ecole Nationale Sup\'erieure Polytechnique,\\
B.P.~8390 Yaound\'e, Cameroun\\
e-mail:tbouetou@polytech.uninet.cm}

\maketitle \begin{quote}{\bf Abstract.}
\em We suggest a new, alternative algebraic method for computation the quantities  $\overset{1}{\nabla_{l}}a_{jk}^{i}$, $\overset{2}{\nabla_{l}}a_{jk}^{i}$ and $d_{jklm}^{i}$ by means of the embedding of local loops into Lie groups.
\end{quote}

\maketitle \begin{quote}{\bf Keywords:}
\em Homogeneous smooth  Loops, Lie groups, Lie algebras, 3-Webs, Chern connection
\end{quote}
\maketitle \begin{quote}{\bf AMS subject classification 2000:}
 11E57, 14C21, 19L99, 20N05, 22E60, 22E67, 32C22, 53A60
\end{quote}

\begin{center}
{\bf Introduction}
\end{center}
  The development of geometry of fiber bundles and foliations stimulates the
interest for new investigation of three-Webs. In \cite{ak3,ak5,go}  
 the techniques was developed for webs using the intresic geometry structure.
In this investigation, we propose to give another approach of computation of
some classical relations, using the technique of the projective space. Our
approach is based on the embedding of a smooth loop into a Lie group, by means of a closed subgroup. This transports the geometric problem into an abstract 
algebraic problem, where the 3-Web is seen as a homogeneous space coset in a generic position. Using this technique the computation of the tensor structure of local loop yield. Therefore we give an application of the computation of the well known tensor 
\em We use algebraic methods to compute the relations $\overset{1}{\nabla_{l}}a_{jk}^{i}$, $\overset{2}{\nabla_{l}}a_{jk}^{i}$ and $d_{jklm}^{i}.$

\section{Analytic representation of law of composition of local smooth loops, embedding in Lie groups}

Let $<G, \cdot, e>$ be a local Lie group and let $H$ be its local closed subgroup.
 Denote by $\mathfrak{G}$ and $\mathfrak{h}$ their corresponding Lie algebra and
 Lie subalgebra and let $Q$ be a smooth space section of left coset $G mod H$ 
passing through $e$ the unity element of $G (e \in G)$.\\
The composition law:\\
$$ \times :Q\times Q \longrightarrow Q$$
$$ (x,y) \longmapsto x\times y=\prod_{Q}(x \cdot y),$$
where $\prod_{Q}: G \to Q$ is the projection on $Q$ parallel to the subgroup $H$, defines in $Q$ a structure of a local loop, i.e  $<Q, \times ,e>$-loop  \cite{Ba,sam1,sab1} .\\
Let us map the tangent space $T_{e}Q$ with the vector subspace $V \subset G$. Then $\mathfrak{G}
=V \dotplus \mathfrak{h}$ since the submanifolds $Q$ and $H$ are transversal in the Lie group $G$.\\
Let us introduce the mapping $\phi$:\\
$$\phi : V  \longrightarrow \mathfrak{h}$$
$$ \xi \longmapsto \phi (\xi)$$ 
defined by the condition $exp (\xi +\phi(\xi)) \in Q$ (for every vector $\xi \in V$, in the neighborhood of $O$, and the map $\phi$ is well defined).\\ Then $\phi(O)=O$ and\\
$$\phi(\xi)=R(\xi , \xi)+S(\xi, \xi, \xi)+o(3)$$
where $$R: V \times V \longrightarrow \mathfrak{h} $$
 $$S: V \times V \times V  \longrightarrow \mathfrak{h} \eqno{(1.1)}$$

are bilinear and trilinear symmetric maps. A base $<e_1, e_2,....,e_N>$ is fixed in $\mathfrak{G}$ such that $<e_1, e_2,....,e_n>$ generate $V$ i.e. $V=<e_1, e_2,....,e_n>$ and $<e_{n+1}, e_{n+2},....,e_N>$ generate $ \mathfrak{h}$:  $\mathfrak{h}=<e_{n+1}, e_{n+2},....,e_N>$. Introduce in the local Lie group $G$ the following normal coordinates, the coordinate on the submanifold $Q$ which is the projection
 from $\exp V$, that is for all $x \in Q$, $x=(x^i)_{i=\overline{1,n}}$, this 
mean $\exp (x^i e_i +\phi(x^i e_i ))=x \in Q$\\
\\
Introduce  the map $$Q \longrightarrow V$$
$$x\longmapsto \overline{x}=x^i e_i. $$ Then the condition written before is equivalent to \\ $$ \overline{x}+\phi ( \overline{x})=x \in Q. $$
In what follows, we will compute the constructed coordinates, fixed on the 
submanifold $Q$.\\
It is known that the law of composition in a Lie group $G(\cdot)$ has  the following representation up to the fourth order in the normal coordinates:\\

$$a \cdot b=a+b+\frac{1}{2}[a,b]+\frac{1}{12}[a,[a,b]]+\frac{1}{12}[b,[b,a]]$$
$$\text{               }-\frac{1}{48}[b,[a,[a,b]]]-\frac{1}{48}[a,[b,[a,b]]]+o(4). \eqno{(1.1)'}$$

Consider the coordinate  representation of the law of composition  $ \times $;for  $y$: $ x=( \overline{x})$ and  $ y=( \overline{y})$ in $Q$. We have:\\

$$
\overline{(x \times y)}= \overline{x}+ \overline{y}+K( \overline{x}, \overline{y})+L( \overline{x}, \overline{x}, \overline{y})+M( \overline{x}, \overline{y}, \overline{y})+$$

 $$\textrm{                                                           } P( \overline{x}, \overline{x}, \overline{x}, \overline{y})+Q( \overline{x}, \overline{x}, \overline{y}, \overline{y})+U( \overline{x}, \overline{y}, \overline{y}, \overline{y})+o(4) \eqno{(1.2)} $$

(Our notation are similar to the notations of  the work \cite{aksh}).\\
Denote the right side in (1.2) by $z=( \overline{z})$. Then for its computation we obtain the equation 
$$\exp ( \overline{z}+\phi ( \overline{z}))=\exp ( \overline{x}+\phi ( \overline{x}))\cdot \exp ( \overline{y}+\phi ( \overline{y}))h  \eqno{(1.3)} $$

where $h$ is and element from $\mathfrak{h}$ in deed we have $h=h( \overline{x}, \overline{y})$.\\
The following proposition holds:\\

{\bf Proposition 1.1}  We have:\\
$$K( \overline{x}, \overline{y})=\frac{1}{2}\prod [ \overline{x}, \overline{y}]$$
\\
where $\prod [ \overline{x}, \overline{y}]$ is the projection of the commutator $ [ \overline{x}, \overline{y}]$ on $V$ parallel to the subalgebra $\mathfrak{h}$.\\
$$h(x,y)= -\frac{1}{2}[ \overline{x}, \overline{y}]+\frac{1}{2}\prod [ \overline{x}, \overline{y}]+2R(\overline{x}, \overline{y})+o(2).$$

Proof: we use the formulae (1.3). Comparing the terms from $V$ and $\mathfrak{h}$
and considering only the terms of first order we obtain:
$$ \overline{z}= \overline{x}+ \overline{y} \in V$$
$$h=o \in \mathfrak{h}.$$
For computing the term of second order we denote\\
$$ \overline{z}= \overline{x}+ \overline{y}+K( \overline{x}, \overline{y}) \in V$$
$$h=N( \overline{x}, \overline{y}) \in \mathfrak{h}$$
from (1.3) and considering (1.1) and (1.1)' we have:\\

$$ \overline{x}+ \overline{y}+K( \overline{x}, \overline{y})+R( \overline{x}, \overline{x})+R( \overline{y}, \overline{y})+2R( \overline{x}, \overline{y})=\overline{x}+ \overline{y}+N( \overline{x}, \overline{y})+R( \overline{x}, \overline{x})+R( \overline{y}, \overline{y})+\frac{1}{2}[ \overline{x}, \overline{y}]$$
then by comparing term from $V$ and $\mathfrak{h}$ and noting that:\\

$$\frac{1}{2}[ \overline{x}, \overline{y}]=\frac{1}{2}\prod [ \overline{x}, \overline{y}]+(\frac{1}{2}[ \overline{x}, \overline{y}]-\frac{1}{2}\prod [ \overline{x}, \overline{y}])$$
hence  $$K( \overline{x}, \overline{y})=\frac{1}{2}\prod [ \overline{x}, \overline{y}]$$
\\
$$h(x,y)= -\frac{1}{2}[ \overline{x}, \overline{y}]+\frac{1}{2}\prod [ \overline{x}, \overline{y}]+2R(\overline{x}, \overline{y})$$

{\bf Corollary 1.1:} from the proposition (1.1) it follows that $$ \overline{(x \times y)}= \overline{x}+ \overline{y}+\frac{1}{2}\prod [ \overline{x}, \overline{y}]+o(2)$$

{\bf Proposition 1.2}  One can show:\\
 $$L( \overline{x}, \overline{x}, \overline{y})=-\frac{1}{6}\prod [ \overline{x}, [ \overline{x}, \overline{y}]]+\frac{1}{2}\prod [R( \overline{x}, \overline{x}), \overline{y}]+\ \frac{1}{4}\prod [ \overline{x}, \prod [ \overline{x}, \overline{y}]]+ \prod [  \overline{x},R( \overline{x}, \overline{y})]$$
\\
\\
$$ M( \overline{x}, \overline{y}, \overline{y})=\frac{1}{3}\prod [ \overline{y},[ \overline{y}, \overline{x}]]+\frac{1}{2}\prod [ \overline{x},R( \overline{y}, \overline{y})]-\frac{1}{4}\prod [ \overline{y},\prod [ \overline{y}, \overline{x}]]+\prod [ \overline{y},R( \overline{x}, \overline{y})]$$
\\
\\
$$ h( \overline{x}, \overline{y})=-\frac{1}{2}[\overline{x}, \overline{y}]+\frac{1}{2} \prod [\overline{x}, \overline{y}]+2R(\overline{x}, \overline{y})+R(\overline{x},\prod [\overline{x}, \overline{y}])+3S(\overline{x},\overline{x},\overline{y})+$$

$$\frac{1}{6}\Lambda [\overline{x},[\overline{x},\overline{y}]]-\frac{1}{4}\Lambda [\overline{x}, \prod [\overline{x},\overline{y}]]-\frac{1}{2} \Lambda [R(\overline{x},\overline{x}),\overline{y}]- \Lambda [\overline{x},R(\overline{x},\overline{y})]+$$

$$+R(\overline{y},\prod [\overline{x}, \overline{y}])+3S(\overline{x},\overline{y},\overline{y})-\frac{1}{3}\Lambda [\overline{y},[\overline{y},\overline{x}]]+$$
$$\frac{1}{4}\Lambda [\overline{y},\prod [\overline{y},\overline{x}]]-\frac{1}{2}\Lambda [\overline{x},R(\overline{y},\overline{y})]-\Lambda [\overline{y},R(\overline{x},\overline{y})]+0(3)$$

where $\Lambda: \mathfrak{G} \longrightarrow \mathfrak{h}$ is the projection on $\mathfrak{h}$ parallel to $V$.\\

Proof. The proof is based on the direct computation. 

Denote:\\

$$\overline{z}= \overline{x}+ \overline{y}+\frac{1}{2}[ \overline{x}, \overline{y}]+L( \overline{x}, \overline{x}, \overline{y})+M( \overline{x}, \overline{y} , \overline{y})$$
and
$$ h(\overline{x},\overline{y})=-\frac{1}{2}[ \overline{x}, \overline{y}]+\frac{1}{2} \prod [\overline{x},\overline{y}]+2R( \overline{x}, \overline{y})+E( \overline{x}, \overline{x}, \overline{y})+F( \overline{x}, \overline{y} , \overline{y}).$$

From (1.3) with the consideration of (1.1) and (1.1)' we obtain the equation\\
$$L(\overline{x}, \overline{x}, \overline{y})+M(\overline{x}, \overline{y}, \overline{y})+R(\overline{x}, \prod [\overline{x}, \overline{y}])+R(\overline{y}, \prod [\overline{x}, \overline{y}])+S(\overline{x}, \overline{x}, \overline{x})+3S(\overline{x}, \overline{y}, \overline{y})+$$ 
$$3S(\overline{x}, \overline{x}, \overline{y})+S(\overline{y}, \overline{y}, \overline{y})+.....=\frac{1}{12} [\overline{x}, [\overline{x}, \overline{y}]]+\frac{1}{12} [\overline{y}, [\overline{y}, \overline{x}]]+E(\overline{x}, \overline{x}, \overline{y})+F(\overline{x}, \overline{y}, \overline{y})+$$ 
$$S(\overline{x}, \overline{x}, \overline{x})+S(\overline{y}, \overline{y}, \overline{y})+\frac{1}{2}[R(\overline{x}, \overline{x}), \overline{y}]+\frac{1}{2} [\overline{x},R( \overline{y}, \overline{y})]+\frac{1}{4}[\overline{x}+ \overline{y}, \prod [\overline{x}, \overline{y}]]-$$ 
$$\frac{1}{4}[\overline{x}+ \overline{y}, [\overline{x}, \overline{y}]]+[\overline{x}+ \overline{y}, R(\overline{x}, \overline{y})]+....$$

Then by comparing term from $V$ and $\mathfrak{h}$ in the last identity we obtain the requirement for $L(\overline{x},\overline{x}, \overline{y})$, $M(\overline{x},\overline{y}, \overline{y})$ and $h(\overline{x}, \overline{y})$ 

in addition 

$$E( \overline{x}, \overline{x} , \overline{y})=R( \overline{x}, \prod [\overline{x} , \overline{y}])+3S( \overline{x}, \overline{x} , \overline{y})+\frac{1}{6} \Lambda [ \overline{x},[ \overline{x} , \overline{y}]]-\frac{1}{4} \Lambda [ \overline{x},\prod [ \overline{x} , \overline{y}]]-$$

$$-\frac{1}{2} \Lambda [ R(\overline{x}, \overline{x}) , \overline{y}]- \Lambda [ \overline{x},R( \overline{x} , \overline{y})] \eqno{(1.4)}  $$

$$F( \overline{x}, \overline{y} , \overline{y})=R( \overline{y}, \prod [\overline{x} , \overline{y}])+3S( \overline{x}, \overline{y} , \overline{y})-\frac{1}{3} \Lambda [ \overline{y},[ \overline{y} , \overline{x}]]+\frac{1}{4} \Lambda [ \overline{y},\prod [ \overline{y} , \overline{x}]]-$$

$$-\frac{1}{2} \Lambda [ \overline{x}, R(\overline{y} , \overline{y})]- \Lambda [ \overline{y},R( \overline{x} , \overline{y})] \eqno{(1.5)} $$
\\
\\
{{\bf Corollary 1.2:} One can obtain:

$$ \overline{(x \times y)}= \overline{x}+ \overline{y}+\frac{1}{2}\prod [ \overline{x}, \overline{y}]-\frac{1}{6}\prod [ \overline{x}, [\overline{x} , \overline{y}]]+\frac{1}{2}\prod [R( \overline{x}, \overline{x}) , \overline{y}]+\frac{1}{4}\prod [ \overline{x}, \prod [\overline{x} , \overline{y}]]+$$

$$+\prod [ \overline{x}, R(\overline{x} , \overline{y})] +\frac{1}{3}\prod [ \overline{y}, [\overline{y} , \overline{x}]]+\frac{1}{2}\prod [ \overline{x},R(\overline{y} , \overline{y})]- $$

$$-\frac{1}{4}\prod [ \overline{y}, \prod [\overline{y} , \overline{x}]]+ \prod [ \overline{y}, R(\overline{x} , \overline{y})]+ o(3).  \eqno{(1.6)}$$
\\  

For the computation of terms of fourth order, denote\\
\\

$$\overline{z}=(1.6)+P( \overline{x}, \overline{x} ,\overline{x}, \overline{y})+Q( \overline{x}, \overline{x} ,\overline{y}, \overline{y})+U( \overline{x}, \overline{y} ,\overline{y}, \overline{y})$$
\\

and for $h$ to take terms of third order.
\\

$P( \overline{x}, \overline{x} ,\overline{x}, \overline{y})+Q( \overline{x},  \overline{x} ,\overline{y}, \overline{y})+U( \overline{x}, \overline{y} ,\overline{y}, \overline{y})=[\overline{x}+R(\overline{x},\overline{x})+S(\overline{x},\overline{x},\overline{x})]\cdot [\overline{y}+R(\overline{y},\overline{y})+S(\overline{y},\overline{y},\overline{y})]\cdot (-\frac{1}{2}\Lambda [\overline{x},\overline{y}]+2R(\overline{x},\overline{y})+E(\overline{x},\overline{x},\overline{y})+F(\overline{x},\overline{y},\overline{y})+....) $
\\

in the fourth order one needs to compute only the term in $V$. Conducting the reasoning as in the past cases one obtain:\\

$P( \overline{x}, \overline{x} ,\overline{x}, \overline{y})+Q( \overline{x},  \overline{x} ,\overline{y}, \overline{y})+U( \overline{x}, \overline{y} ,\overline{y}, \overline{y})=\Bigg\{[\overline{x}+R(\overline{x},\overline{x})+S(\overline{x},\overline{x},\overline{x})+\overline{y}+R(\overline{y},\overline{y})+S(\overline{y},\overline{y},\overline{y})+ \frac{1}{2}[\overline{x},\overline{y}]+\\+\frac{1}{2}[\overline{x},R(\overline{y},\overline{y})]+\frac{1}{2}[R(\overline{x},\overline{x}),R(\overline{y},\overline{y})]+\frac{1}{2}[\overline{x},S(\overline{y},\overline{y},\overline{y})]+\frac{1}{2}[S(\overline{x},\overline{x},\overline{x}),\overline{y})]+\frac{1}{12}[\overline{x},[\overline{x},\overline{y}]]+\\+\frac{1}{12}[\overline{x},[\overline{x},R(\overline{y},\overline{y})]]+\frac{1}{12}[\overline{y},[\overline{y},\overline{x}]]+ \frac{1}{12}[\overline{y},[\overline{y},R(\overline{x},\overline{x})]]-\frac{1}{48}[\overline{y},[\overline{x},[\overline{x},\overline{y}]]]- \\-\frac{1}{48}[\overline{x},[\overline{y},[\overline{x},\overline{y}]]]+...\Bigg\}\cdot (-\frac{1}{2}\Lambda [\overline{x},\overline{y}]+2R(\overline{x}.\overline{y})+E(\overline{x},\overline{x},\overline{y})+F(\overline{x},\overline{y},\overline{y})+...)=   $

$$= \frac{1}{2}\prod [\overline{x},E(\overline{x},\overline{x},\overline{y})]+ \frac{1}{2}\prod [\overline{x},F(\overline{x},\overline{y},\overline{y})]+ \frac{1}{2}\prod [\overline{y},E(\overline{x},\overline{x},\overline{y})]+ \frac{1}{2}\prod [\overline{y},F(\overline{x},\overline{y},\overline{y})]-$$ 
$$\frac{1}{8}\prod [\prod [\overline{x},\overline{y}],[\overline{x},\overline{y}]]+ \frac{1}{2}\prod[\prod [\overline{x},\overline{y}],R(\overline{x},\overline{y})]+\frac{1}{12}\prod [\overline{x},[\overline{x},-\frac{1}{2}\Lambda [\overline{x},\overline{y}]+2R(\overline{x},\overline{y})]]$$ 
$$+\frac{1}{12}\prod [\overline{y},[\overline{y},-\frac{1}{2}\Lambda [\overline{x},\overline{y}]+2R(\overline{x},\overline{y})]]+\frac{1}{12}\prod [\overline{x},[\overline{y},-\frac{1}{2}\Lambda [\overline{x},\overline{y}]+2R(\overline{x},\overline{y})]]+$$ 
$$\frac{1}{12}\prod [\overline{y},[\overline{x},-\frac{1}{2}\Lambda [\overline{x},\overline{y}]+2R(\overline{x},\overline{y})]]+\frac{1}{2}\prod [\overline{x},S(\overline{y},\overline{y},\overline{y})]+\frac{1}{2}\prod [S(\overline{x},\overline{x},\overline{x}),\overline{y}]+$$ 
$$\frac{1}{12}\prod [\overline{x},[\overline{x},R(\overline{y},\overline{y})]]+\frac{1}{12}\prod [\overline{y},[\overline{y},R(\overline{x},\overline{x})]]-\\-\frac{1}{48}\prod [\overline{y},[\overline{x},[\overline{x},\overline{y}]]]-\frac{1}{48}\prod [\overline{x},[\overline{y},[\overline{x},\overline{y}]]].$$
\\
all the equality in the above expression are modulo $\mathfrak{h}$.
\\

Then the following proposition holds:
\\
\\
{\bf Proposition 1.3}\\
$$P(\overline{x},\overline{x},\overline{x},\overline{y})=-\frac{1}{2}\prod [\overline{y},S(\overline{x},\overline{x},\overline{x})]+\frac{1}{12}\prod [\overline{x},[\overline{x},-\frac{1}{12}\Lambda [\overline{x},\overline{y}]+2R(\overline{x},\overline{y})]]+$$

$$\frac{1}{2}\prod [\overline{x},E(\overline{x},\overline{x},\overline{y})]  \eqno{(1.7)}$$

$$U(\overline{x},\overline{y},\overline{y},\overline{y})=\frac{1}{2}\prod [\overline{x},S(\overline{y},\overline{y},\overline{y})]+\frac{1}{12}\prod [\overline{y},[\overline{y},-\frac{1}{12}\Lambda [\overline{x},\overline{y}]+2R(\overline{x},\overline{y})]]+$$

$$+\frac{1}{2}\prod [\overline{y},F(\overline{x},\overline{y},\overline{y})]  \eqno{(1.8)}$$

$$Q(\overline{x},\overline{x},\overline{y},\overline{y})=\frac{1}{2}\prod [\overline{y},E(\overline{x},\overline{x},\overline{y})] +\frac{1}{2}\prod [\overline{x},F(\overline{x},\overline{y},\overline{y})]-\frac{1}{8}\prod [\prod [\overline{x},\overline{y}],[\overline{x},\overline{y}]]+$$

$$+\frac{1}{2}\prod [\prod [\overline{x},\overline{y}],R(\overline{x},\overline{y})]+ \frac{1}{12}\prod [\overline{x},[\overline{y},-\frac{1}{2}\Lambda [\overline{x},\overline{y}]+2R(\overline{x},\overline{y})]]+$$

$$+\frac{1}{12}\prod [\overline{y},[\overline{x},-\frac{1}{2}\Lambda [\overline{x},\overline{y}]+2R(\overline{x},\overline{y})]]+\frac{1}{12}\prod [\overline{x},[\overline{x},R(\overline{y},\overline{y})]]+$$

$$+\frac{1}{12}\prod [\overline{y},[\overline{y},R(\overline{x},\overline{x})]]
-\frac{1}{48}\prod [\overline{y},[\overline{x},[\overline{x},\overline{y}]]]-\frac{1}{48}\prod [\overline{x},[\overline{y},[\overline{x},\overline{y}]]] \eqno{(1.9)} $$
\\
\\
{\bf Corollary 1.3:}\\

$$
( \overline{x \times y} )=\overline{x}+\overline{y}+\frac{1}{2}\prod [\overline{x},\overline{y}]-\frac{1}{6}\prod [\overline{x},[\overline{x},\overline{y}]]+\frac{1}{2}\prod [R(\overline{x},\overline{x}),\overline{y}]+\frac{1}{4}\prod [\overline{x}, \prod [\overline{x},\overline{y}]]+$$

$$+\prod [\overline{x},R(\overline{x},\overline{y})]+\frac{1}{3}\prod [\overline{y},[\overline{y},\overline{x}]]+\frac{1}{2}\prod [\overline{x},R(\overline{y},\overline{y})]-\frac{1}{4}\prod [\overline{y}, \prod [\overline{y},\overline{x}]]+$$

$$+\prod [\overline{y},R(\overline{x},\overline{y})]+P( \overline{x}, \overline{x} ,\overline{x}, \overline{y})+Q( \overline{x}, \overline{x} ,\overline{y}, \overline{y})+U( \overline{x}, \overline{y} ,\overline{y}, \overline{y})+0(4)  \eqno{(1.10)}$$

Where $P( \overline{x}, \overline{x} ,\overline{x}, \overline{y}),Q( \overline{x}, \overline{x} ,\overline{y}, \overline{y})$  and $U( \overline{x}, \overline{y} ,\overline{y}, \overline{y})$ are from (1.7), (1.8) and (1.9)

\section{Tensor structure of a smooth analytic loop}

Let $<Q, \times, e>$ be a smooth analytic loop with the neutral element $e$. In  a standard way see \cite{Be} on the Cartesian product $Q \times Q$ we introduce the structure of a three-webs $W$ such that the  submanifold in the view of $\{a\} \times Q$ is a vertical foliations $(a \in Q)$, $Q \times \{b\}$ is a horizontal foliations $(b \in Q)$ and  the set $ \{ (a,b): a \times b =c=conts \}$ the foliations of the third family $(c \in Q)$. In the coordinate $(x^1, x^2,.....,x^n, y^1, y^2, .......,y^n )$,  the indicated foliations are described by the system of differential 1-form \cite{ak, Mih}.\\

$$ \omega^i_1=o,  \omega^i_2=o,  \omega^i_3=  \omega^i_1+ \omega^i_2=o$$
where $$ \omega^i_1= P^i_{\alpha} dx^{\alpha}, \;  \omega^i_{2}= Q^i_{\beta} dy^{\beta},$$ $$P^i_{\alpha} (x,y)=\frac{\partial \mu^i}{\partial x^{\alpha}}$$
 $$Q^i_{\beta} (x,y)=\frac{\partial \mu^i}{\partial y^{\beta}}$$
$$\mu^i (x,y)=(x \times y)^i$$

In the space of a 3-Web $W$, introduce the so called Chern canonical connection\\ $\nabla=(\overset{1}{\nabla}, \overset{2}{\nabla})$ \cite{ aksh,she2}.

The indicated connection is described by:

$$ \omega^{k}_{j}=\Gamma^{k}_{ij} \omega^{i}_{1}+\Gamma^{k}_{jl} \omega^{j}_{2}, $$
$$\Gamma^{k}_{ij}=-\widetilde{P}^{\alpha}_{i}\widetilde{Q}^{\beta}_{j}\frac{\partial^{2}\mu^{k}}{\partial x^{\alpha} \partial y^{\beta}}$$

where $\widetilde{P}^{\alpha}_{i} $ and $\widetilde{Q}^{\beta}_{j} $ are inverse matrices for  $P^{\alpha}_{i} $ and $Q^{\beta}_{j} $ respectively in  terms of the following structural equations:

$$ d\omega^{k}_{1}= \omega^{l}_{1}\wedge \omega^{k}_{l}+a^{k}_{ij} \omega^{i}_{1}\wedge \omega^{j}_{l}$$

$$ d\omega^{k}_{2}= \omega^{l}_{2}\wedge \omega^{k}_{l}-a^{k}_{ij} \omega^{i}_{2}\wedge \omega^{j}_{2} \eqno{(2.1)}$$

$$ d\omega^{k}_{j}= \omega^{i}_{j}\wedge \omega^{k}_{i}+b^{k}_{jlm} \omega^{l}_{1}\wedge \omega^{m}_{2}$$   

where $$ a^{k}_{ij}=-\frac{1}{2}\frac{\partial^{2}\mu^{k}}{\partial x^{\alpha} \partial y^{\beta}}(\widetilde{P}^{\alpha}_{i}\widetilde{Q}^{\beta}_{j}-\widetilde{P}^{\alpha}_{j}\widetilde{Q}^{\beta}_{i}) $$

$$b^{k}_{jlm}=(-\frac{\partial^{3} \mu^{k}}{\partial x^{\alpha}\partial x^{\beta}\partial y^{\gamma}}\widetilde{P}^{\beta}_{j}+\frac{\partial^{3} \mu^{k}}{\partial x^{\alpha}\partial y^{\beta}\partial y^{\gamma}}\widetilde{Q}^{\beta}_{j})\widetilde{P}^{\alpha}_{l}\widetilde{Q}^{\gamma}_{m}-\Gamma^{k}_{pm}\frac{\partial^{2}\mu^{p}}{\partial x^{\alpha} \partial x^{\beta}}\widetilde{P}^{\alpha}_{l}\widetilde{P}^{\beta}_{j}+$$

$$+\Gamma^{k}_{lp}\frac{\partial^{2}\mu^{p}}{\partial y^{\alpha} \partial y^{\beta}}\widetilde{P}^{\alpha}_{j}\widetilde{Q}^{\beta}_{m}-\Gamma^{k}_{pm}\Gamma^{p}_{lj}+\Gamma^{k}_{lp}\Gamma^{p}_{jm}$$

The Chern connection in the 3-Web associated to the loop $<Q, \times ,e>$, 
admits an
 alternative description in  terms of anti-product of the loop $Q$ by itself \cite{sam3,sab2}. In the set $Q\times Q$ introduce the covering loopouscular structure, by denoting for any pair $X=(x,x')$, $Y=(y,y')$, $A(u,v)$

$$L(X,A,Y)=((x(u\backslash yv))/v,u\backslash ((uy'/v)x')) \eqno{(2.2)}.$$

Then the Chern connection coincide with the connection tangent to the covering loopuscular structure \cite{sab2}.

In particular, for any tensor field $\Omega (u,v)$, in the space of 3-web $W=Q\times Q$

$$\overset{1}{\nabla_{i}}\Omega (u=e, v=e)=\frac{\partial}{\partial u^i}\Bigg [ \Big \{ [L^{(e,e)}_{(u,e)}]_{\ast,(e,e)}\Big \}^{-1}\Omega (u,e)\Bigg ]\arrowvert_{ u=e} \eqno{(2.3)},$$ 

$$\overset{2}{\nabla_{i}}\Omega (u=e, v=e)=\frac{\partial}{\partial v^i}\Bigg [ \Big \{ [L^{(e,e)}_{(e,v)}]_{\ast,(e,e)}\Big \}^{-1}\Omega (e,v)\Bigg ]\arrowvert_{ v=e} .$$

The value in the point $(e,e)$ of the 3-Web $W= Q \times Q$  to the loop
$< Q, \times, e>$ the fundamental tensor field $a^{i}_{jk}$, $b^{i}_{jkl}$ and
their corresponding derivations $ \overset{1}{\nabla_{i}}$, 
$\overset{2}{\nabla_{i}}$ are called  the structure tensors of the  loop. The 
structure tensor of the smooth loop $ <Q, \times, e  >$
defined uniquely by  its construction up to isomorphism \cite{aksh,Mih,sam1,she2}\\
\\
\\
{\bf Proposition 2.1} \cite{ak, she2} The following relations hold:\\
$$\overset{1}{\nabla_{l}}a^{i}_{jk}=b^{i}_{[j\arrowvert l\arrowvert k]} \eqno {(2.4)}$$

$$\overset{2}{\nabla_{l}}a^{i}_{jk}=b^{i}_{[jk]l} \eqno {(2.5)}.$$\\

For the proof of the proposition, it's sufficient to consider the first differential expression of the system $(2.1)$.\\

Introduce the notation\\
$$ c^{i}_{jklm}=\overset{1}{\nabla}_{m} b^{i}_{jkl}\arrowvert_{(e,e)} $$
 
$$ d^{i}_{jklm}=\overset{2}{\nabla_{m}}b^{i}_{jkl}\arrowvert_{(e,e)}. $$

And consider the proposition $(1.2)$. The law of composition $(\times)$ of the
smooth local loop  $ <Q, \times, e  >$ in the coordinate $x=(\overline{x})$ 
centralized at the point $e$, is given by:\\

$$\overline{(x \times y)}=\overline{x}+\overline{y}+K(\overline{x},\overline{y})+L(\overline{x},\overline{x},\overline{y})+M(\overline{x},\overline{y},\overline{y})+P(\overline{x},\overline{x},\overline{x},\overline{y})$$

$$+Q(\overline{x},\overline{x},\overline{y},\overline{y})+U(\overline{x},\overline{y},\overline{y},\overline{y})+o(4).$$
\\

Consider  $ <Q, \times, e  >$ as a coordinate loop of the 3-Web $W$, defined in the neighbourhood of the point $(e,e)$ of the manifold $Q \times Q$. Then in 
conformity with \cite{aksh, she1} the basic tensor of the web can be 
expressed in term the of coefficient of the decomposition of the loop in the following way:\\

$$ a(\overline{x},\overline{y})=-K(\overline{x},\overline{y}),$$
$$ b(\overline{x},\overline{y},\overline{z})=-B(\overline{y},\overline{x},\overline{z})  \eqno {(2.6)}$$
\\

$$ c(\overline{x},\overline{y},\overline{z},\overline{t} )=(4Q-6P)(\overline{y},\overline{t},\overline{x},\overline{z})+a(\overline{t},b(\overline{x},\overline{y},\overline{z}))+a(\overline{y},b(\overline{x},\overline{t},\overline{z}))$$

$$-b(\overline{x},a(\overline{t},\overline{y}),\overline{z})+a(2L(\overline{y},\overline{t},\overline{x}),\overline{z})-2L(a(\overline{x},\overline{y}),\overline{t},\overline{z})$$

$$-2L(\overline{y},a(\overline{x},\overline{t}),\overline{z})-2L(\overline{y},\overline{t},a(\overline{x},\overline{z})) \eqno{(2.7)}$$\\

$$d(\overline{x},\overline{y},\overline{z},\overline{t})=(4Q-6P)(\overline{y},\overline{x},\overline{z},\overline{t})-a(b(\overline{x},\overline{y},\overline{z}),\overline{t})-a(b(\overline{x},\overline{y},\overline{t}),\overline{z})+$$

$$+b(\overline{x},\overline{y},a(\overline{z},\overline{t}))+a(\overline{y},2M(\overline{x},\overline{z},\overline{t}))-2M(a(\overline{y},\overline{x}),\overline{z},\overline{t})-$$

$$-2M(\overline{y},a(\overline{z},\overline{x}),\overline{t})-2M(\overline{y},\overline{z},a(\overline{t},\overline{x})) \eqno{(2.8)}$$
\\

where $$B(\overline{x},\overline{y},\overline{z})=2L(\overline{x},\overline{y},\overline{z})-2M(\overline{x},\overline{y},\overline{z})-K(\overline{x},K(\overline{y},\overline{z}))+K(K(\overline{x},\overline{y}),\overline{z}) \eqno{(2.10)}$$
 
\section{Structure tensor of a smooth local loop, Embedding in Lie group}

Let $<Q, \times, e>$ be a local smooth loop, the embedding in the Lie group $G$ as a section of left coset $G mod H$, where $H$ is a closed subgroup in $G$. In
what follows, we will consider that  $<Q, \times, e>$, is referred to the normal coordinates $X=(\overline{x})$.\\
\\

{\bf Proposition 3.1} The following relations holds:\\

$$ a(\overline{x},\overline{y})=-\frac{1}{2}\prod [\overline{x},\overline{y}] \eqno {(3.1)}$$
\\
$$b(\overline{x},\overline{y},\overline{z})=-\frac{1}{2}\prod [[\overline{x},\overline{y}],\overline{z}]+\frac{1}{2}\prod [\prod [\overline{x},\overline{y}],\overline{z}]-2\prod [R(\overline{x},\overline{y}),\overline{z}] \eqno {(3.2)}$$
\\

Proof: The first relation follows from the proposition 1.1 and the relation $(2.6)$. In the relation $(2.10)$ we have:

$$B(\overline{x},\overline{y},\overline{z})=2L(\overline{x},\overline{y},\overline{z})-2M(\overline{x},\overline{y},\overline{z})-K(\overline{x},K(\overline{y},\overline{z}))+K(K(\overline{x},\overline{y}),\overline{z}) $$
\\
and from the proposition 1.2 we have:\\

$$2L(\overline{x},\overline{y},\overline{z})=-\frac{1}{6}\prod [\overline{x},[\overline{y},\overline{z}]]+\prod[R(\overline{x},\overline{y}),\overline{z}]+\frac{1}{4}\prod [\overline{x},\prod[\overline{y},\overline{z}]]-\frac{1}{6}\prod [\overline{y},[\overline{x},\overline{y}]]+$$

$$+\prod [\overline{x},R(\overline{y},\overline{z})]+\prod [\overline{y},R(\overline{x},\overline{z})]+\frac{1}{4}\prod [\overline{y},\prod [\overline{x},\overline{z}]]$$\\

$$2M(\overline{x},\overline{y},\overline{z})=\frac{1}{3}\prod [\overline{y},[\overline{z},\overline{x}]]+\prod[\overline{x},R(\overline{y},\overline{z})]-\frac{1}{4}\prod [\overline{y},\prod [\overline{z},\overline{x}]]+\frac{1}{3}\prod [\overline{z},[\overline{y},\overline{x}]]+$$

$$+\prod [\overline{y},R(\overline{x},\overline{z})]+\prod [\overline{z},R(\overline{x},\overline{y})]-\frac{1}{4}\prod  [\overline{z},\prod [\overline{y},\overline{x}]]$$
\\

further more $$K(\overline{x},K(\overline{y},\overline{z}))=\frac{1}{4}\prod  [\overline{x},\prod [\overline{y},\overline{z}]].$$

 $$K(K(\overline{x},\overline{y}),\overline{z}))=\frac{1}{4}\prod  [\prod [\overline{x},\overline{y}],\overline{z}].$$

Substituting these expressions in $B(\overline{x},\overline{y},\overline{z})$, we obtain:\\
$$B(\overline{x},\overline{y},\overline{z})=-\frac{1}{2}\prod [[\overline{x},\overline{y}],\overline{z}]+\frac{1}{2}\prod [\prod [\overline{x},\overline{y}],\overline{z}]+2\prod [R(\overline{x},\overline{y}),\overline{z}]$$

but from $(2.6)$ we have $b(\overline{x},\overline{y},\overline{z})=-B(\overline{y},\overline{x},\overline{z})$.\\
Hence :\\

$$b(\overline{x},\overline{y},\overline{z})=-\frac{1}{2}\prod [[\overline{x},\overline{y}],\overline{z}]-\frac{1}{2}\prod [\prod [\overline{x},\overline{y}],\overline{z}]-2\prod [R(\overline{x},\overline{y}),\overline{z}].$$
\\

Let $\Omega$ be one of the structural tensor of the loop $Q$, and consider the 
expression of the fundamental tensor field $\Omega(u,v)$ in the space of 
three-webs $W=Q\times Q$. Then $\Omega =\Omega (u=e,v=e)$ and for 
$\overset{1}{\nabla}_{i}\Omega (u=e,v=e)$,$ \overset{2}{\nabla}_{i}\Omega (u=e,v=e) $ the formulae obtained in $(2.3)$ hold.\\
Consider the computation of $ \overset{1}{\nabla}_{i}\Omega (u=e,v=e)$,
the value of the tensor field $\Omega (u,v)$ for $v=e$ can be seen as the 
structure of the smooth local loop $<Q, \underset{u}{\times}, u>$ where
$$x  \underset{u}{\times}  y=x\times (u\backslash y).$$ As a result, $\nabla$ is transported from $T_{u}Q$
in  $T_{e}Q$ by means of the inverse transformation $R_u$, which coincide with
the structure of the tensor $\widetilde{\Omega_u}$ and the smooth local loop
 $<Q, \underset{u}{\cdot}, e>$ with the operation:\\

$$x  \underset{u}{\cdot}  y=u\backslash ((u \times x) \times y). \eqno {(3.3)}$$\\
So that  $$ \overset{1}{\nabla}_{i}\Omega (u=e,v=e)=\frac{\partial \widetilde {\Omega_{u}}}{\partial u^{i}}\arrowvert _{u=e}$$ in addition the law of 
composition $(3.3)$ allow an intuitive algebraic interpretation in terms of
the enveloping Lie group $G$.

Consider the section $Q'_{u}=Q \cdot u^{-1}$ of the coset space 
$G/\widetilde{H_{u}}$ where \mbox{ $\widetilde{H_{u}}=u\cdot H\cdot u^{-1}$},
 $u\in Q $ and the 
map: 
$$\Psi_{u} : Q  \longrightarrow Q'_{u}$$
$$ x \longmapsto (u\times x)\times u^{-1}.$$
 Denote by $(  \underset{u}{\ast}  )$ the law of composition in $Q'_{u}$, so that:
$$a \underset{u}{\ast}b=\prod'_{u}(ab) $$ where 

$\prod_{u}':G  \longrightarrow Q'_{u} $ is the projection on $ Q'_{u}$
parallel to  $ \widetilde{H_{u}}$. The following proposition hold.\\
{\bf Proposition 3.2} The map $\Psi_{u} : Q  \longrightarrow Q'_{u}$ is an isomorphism
of the smooth loops  $<Q, \underset{u}{\cdot}, e>$ and $ <Q_{u}', \underset{u}{\ast},e>$

Proof:\\Let $a=\Psi_{u} x$, $b=\Psi_{u} y$ and $a \underset{u}{\ast}b=\Psi_{u} z$\\
where $x,y,z \in Q$. 

Then
 
 $$a \underset{u}{\ast}b=\prod'_{u}(ab)=\prod'_{u}((u\times x)\cdot u^{-1}\cdot (u \times y)\cdot u^{-1})$$
$(a \underset{u}{\ast}b)\times u\cdot h\cdot u^{-1}=(u\times x)u^{-1}\cdot(u\times y)\cdot u^{-1}.$

Multiplying by $u$ obtain:

$$(a \underset{u}{\ast}b)\times u\cdot h=(u\times x)\times y.$$
Applying the projection to the last equality, we obtain

$$(a \underset{u}{\ast}b)\times u=(u\times x)\times y.$$

Furthermore

$$(a \underset{u}{\ast}b)\times u=(\Psi_{u} z)\times u=(u\times z)\cdot u^{-1}\times u=(u\times x)\times y.$$
Then $z=u\backslash (u\times x)\times y$ and

$$(a \underset{u}{\ast}b)=(\Psi_{u}x) \underset{u}{\ast}(\Psi_{u}y )=\Psi_{u}z=   \Psi_{u}\{ u \backslash (u\times x)\times y\}= \Psi_{u} ( x \underset{u}{\cdot}y).$$
Therefore  $\Psi_{u} ( x \underset{u}{\cdot}y)=(\Psi_{u}x)\ast (\Psi_{u}y).$ \\Hence the result.\\
Similarly we establish that:\\

 $$ \overset{2}{\nabla}_{i}\Omega (u=e,v=e)=\frac{\partial \widetilde{\widetilde {\Omega_{v}}}}{\partial v^{i}}\arrowvert _{v=e}$$
where $ \widetilde{\widetilde {\Omega}}$ correspond to the structure tensor of
the local loop $<Q, \frac{1}{v}, e>$ with the composition law:\\
$$ x \frac{1}{v} y=(x\times (y \times v))/v.  \eqno {(3.4)} $$\\

The law of composition $(3.4)$ allows us to find an  algebraic interpretation in
terms of the enveloping Lie group $G$.\\

Let us introduce in consideration the subgroup $H_{v}''=vHv^{-1}$ where
$v \in Q$. The following proposition holds:\\

{\bf Proposition 3.3} $$x \frac{1}{v} y=\prod_{v}''(xy)$$ for all $x,y \in Q$\\ where

$\prod_{v}'':G  \longrightarrow Q $ is the projection on $ Q$
parallel to  $ H_{v}''.$

Proof:\\ In the Lie group $G$ we have $xy=(x\perp y)\times vhv^{-1}$  which is
equivalent to  $xy \cdot v=(x\perp y)\times vh$. Applying $\prod$ to the last
formula we get $$x \times (y\times v)=(x\perp y)\times v.$$ Therefore
$x \perp y=x \times (y \times v)/v$.

\section{Application:Computation of $\stackrel{2}{\nabla}_{l}a^{i}_{jk}$ and $\stackrel{1}{\nabla}_{l}a^{i}_{jk}$}

I:\;  Computation of $ \overset{2}{\nabla}_{l}a^{i}_{jk}$\\

For $u\in Q$, introduce the map

$$Ad_{u} : G  \longrightarrow G$$
$$ x \longmapsto u x u^{-1}.$$
Let $u=\exp \zeta$, where $\zeta \in Q$ and $g \in H$. Then
$$ Ad_{u}(g)=ugu^{-1}=Ad(\exp \zeta)(g)=\exp (ad \zeta (g))$$
$$ = g+[\zeta,g]+o(\zeta)  \eqno {(4.1)}$$
\\
and $ g+[\zeta,g]+o(\zeta)\in H_{u}''$, where $H_{u}''=uHu^{-1}$.\\

Let $\prod''_{u}:\mathfrak{G} \longrightarrow T_{e}Q $ be the projection on
$T_{e}Q $
parallel to $\mathfrak{h}_{u}''$ and $\exp \mathfrak{h}_{u}''=H_{u}''$.\\

By fixing $\xi, \eta $ from $\mathfrak{G}$, we find that\\
$$[\xi, \eta]=\prod [\xi, \eta]+h_{1}  \eqno {(4.2)} $$\\
$$[\xi, \eta]=\prod_{u}'' [\xi, \eta]+h_{2}  \eqno {(4.3)} $$\\
where $h_{1}\in \mathfrak{h}$ and $h_{2}\in \mathfrak{h}_{u}''$. From $(4.1)$ we obtain that $h_{2}$
has the form  $h_{2}= h_{1}+\hat{h}(\zeta)+[\zeta, h_{1}]+o(\zeta)$, where
$\hat{h}(\zeta) \in \mathfrak{h}_{u}''$. From $(4.2)$ and $(4.3)$ it follows that:\\
$$ \prod''_{u}[\xi,\eta ]=[\xi, \eta ]-h_{2}=\prod [\xi, \eta ]-\hat{h}(\zeta)-[\zeta, h_{1}]+o(\zeta)$$
$$=\prod [\xi, \eta ]-\prod [\zeta, h_{1}]+o(\zeta).$$\\

But from $(4.2)$,we have $h_{1}=[\xi, \eta ]-\prod [\xi, \eta ]$. It follows that \\
$$ \prod''_{u}[\xi,\eta ]=\prod [\xi, \eta ]-\prod [\zeta,[\xi, \eta]]+\prod [\zeta, \prod [\xi, \eta ]]+o(\zeta)$$
$$=\prod [\xi,\eta]+\prod [[\xi, \eta],\zeta]-\prod [ \prod [\xi, \eta ], \zeta]+o(\zeta).$$

Denote by $ a_{u}''(\xi, \eta)=-\frac{1}{2}\prod_{u}''[\xi, \eta]$. Then
$$a_{u}''(\xi, \eta)=a(\xi, \eta)-\frac{1}{2}\prod [[\xi, \eta]]+\frac{1}{2}\prod [\prod [\xi,\eta], \zeta].$$

Finally we have:\\
$$\overset{2}{\nabla}_{l}a_{jk}^{i}\xi^{j}\eta^{k}\zeta^{l}=\frac{d}{dt}\Big( a_{\exp  t\zeta}''(\xi, \eta) \Big)\arrowvert_{t=0}=-\frac{1}{2}\prod [[\xi, \eta]]+\frac{1}{2}\prod [\prod [\xi, \eta ], \zeta]$$

We obtain a result in conformity with proposition 2.1 and the relation $(3.2)$
in deed, from the relation $(3.2)$
$$b(\xi, \eta, \zeta)=-\frac{1}{2}[[\xi, \eta], \zeta]+\frac{1}{2}\prod [\prod [\xi, \eta], \zeta]-2\prod [R(\xi, \eta), \zeta ].$$
From which we find
$$\frac{1}{2}[b(\xi, \eta, \zeta)-b(\eta, \xi, \zeta)]=-\frac{1}{2}\prod [[\xi, \eta], \zeta]+\frac{1}{2}\prod [\prod [\xi, \eta], \zeta]$$
so that $\overset{2}{\nabla}_{l}a_{jk}^{i}=b_{[jk]l}^{i}.$\\

   II: \; Computation of $ \overset{1}{\nabla}_{l}a^{i}_{jk}$\\

Let us introduce the map:\\
$$\Psi_{u} : Q  \longrightarrow Q_{u}'$$
$$ x \longmapsto (u \times x) u^{-1}.$$

Then $d\Psi_{u}\arrowvert_{e}  :T_{e}Q  \longrightarrow T_{e}Q_{u}'$.
Then the following proposition holds:\\

{\bf Proposition 4.1} The map define from the tangent space $T_{e}Q$ to tangent space $ T_{e}Q_{u}'$ is defined as follows:\\

 $$d\Psi_{u}\arrowvert_{e}  :T_{e}Q  \longrightarrow T_{e}Q_{u}'$$
$$ \xi \longmapsto \xi+\frac{1}{2}[u,\xi]+\frac{1}{2}\prod [u, \xi]+2R(u,\xi)+o(u).$$

Proof. For the proof of this proposition, using the notion from section 2 and 
the relation $(1.3)$ we have $u \times \xi= (u\cdot \xi)\cdot h$ but from the
proposition 1.4
$$h(u,\xi)=-\frac{1}{2}[u,\xi]+\frac{1}{2}[u,\xi]+2R(u,\xi)+o(u)$$

Thus\\
$$u\times \xi=(u\cdot \xi)\cdot h=u+\xi+\frac{1}{2}\prod [u,\xi]+\frac{1}{2}\prod [u, \xi]+2R(u, \xi)+o(u)$$
and $$ (u\times \xi)\times u^{-1}=u+\xi+\frac{1}{2}\prod [u,\xi]+2R(u,\xi)-u-\frac{1}{2}[\xi,u]+o(u)$$
$$=\xi+\frac{1}{2}\prod [u,\xi]+\frac{1}{2}[u, \xi]+2R(u,\xi)+o(u)$$

Let $\widetilde{\prod_{u}}:\mathfrak{G} \longrightarrow T_{e}Q' $ be the 
projection on $T_{e}Q' $
parallel to $\widetilde{\mathfrak{h}_{u}}'$ where $\exp \widetilde{\mathfrak{h}_{u}}'=uHu^{-1}$.\\

Then we obtain the equation\\
$$\omega+h_{1}=\omega'+h_{1}'+[u,h_{1}']$$ with $\omega \in T_{e}Q$,  
$h_{1}\in \mathfrak{h}$,  $\omega' \in T_{e}Q'$,  $h_{1}'\in \mathfrak{h}$. \\
For the computation of $\omega'=\omega'(u,\omega)$.\\
From the proposition $(4.1)$ we have:
$$\omega +h_{1}=\widetilde{\omega}+\frac{1}{2}\prod [u,\widetilde{\omega}]+\frac{1}{2}[u,\widetilde{\omega}]+2R(u,\widetilde{\omega})+h_{1}'+[u,h_{1}']+o(u)$$
where $\widetilde{\omega} \in T_{e}Q$, so that:
$$\widetilde{\omega}+\frac{1}{2}\prod [u,\widetilde{\omega}]+\frac{1}{2}[u,\widetilde{\omega}]+2R(u,\widetilde{\omega})=\omega'$$
It follows that:
$$\omega=\widetilde{\omega}+\prod [u,\widetilde{\omega}]+[u,h_{1}']$$
$$h_{1}=h_{1}'+\textrm{terms with u} $$
from which 
$$\widetilde{\omega}=\omega-\prod [u, \omega]-[u, h_{1}']$$
$$ h_{1}'=h_{1}+\textrm{term with u} $$
Then substituting in $\omega'$ the expression from $\widetilde{\omega}$ we obtain
that:

$$\omega'=\omega-\prod [u, \omega]-\prod [u,h_{1}]+\frac{1}{2}\prod[u,h_{1}]+\frac{1}{2}[u, \omega]+2R(u, \omega)+o(u)$$
$$=\omega+\frac{1}{2}[u,\omega]-\frac{1}{2}\prod [u, \omega]-\prod [u,h_{1}]+2R(u, \omega)+o(u)$$
from which we find that
$$\widetilde{\prod_{u}}(\omega+h_{1})=\omega'=\omega+\frac{1}{2}[u, \omega]-\frac{1}{2}[u,\omega]+2R(u,\omega)-\prod [u,h_{1}]. \eqno {(4.4)}$$

Now let us compute\\
 $\widetilde{a_{u}}(\xi, \eta)=-\frac{1}{2}(d\Psi)^{-1}\widetilde{\prod_{u}}[d\Psi_{\xi},d\Psi_{\eta}]$

where $\xi, \eta \in T_{e}Q$

$(d\Psi)^{-1}\widetilde{\prod_{u}}[d\Psi_{\xi},d\Psi_{\eta}]=(d\Psi)^{-1}\widetilde{\prod_{u}}\Bigg[\xi+\frac{1}{2}[u,\xi]+\frac{1}{2}\prod [u,\xi]+2R(u, \xi),\eta+\frac{1}{2}[u,\eta]+\frac{1}{2}\prod [u, \eta]+2R(u,\eta)\Bigg] $\\

$=(d\Psi)^{-1}\widetilde{\prod_{u}}\Bigg\{[\xi, \eta]+\frac{1}{2}[\xi,[u,\eta]]+\frac{1}{2}[\xi,\prod [u,\eta]]+2[\xi,R(u,\eta)]-\frac{1}{2}[\eta,[u,\xi]]-\\-\frac{1}{2}[\eta,\prod [u, \xi]]-2[\eta,R(u,\xi)]   \Bigg\}=$\\

$=(d\Psi)^{-1}\Bigg\{\prod [\xi, \eta]+\frac{1}{2}\prod [\xi,[u,\eta]]+\frac{1}{2}\prod [\xi,\prod [u,\eta]]+2\prod [\xi,R(u,\eta)]-\frac{1}{2}\prod [\eta,[u,\xi]]-\frac{1}{2}\prod [\eta,\prod [u, \xi]]-2\prod [\eta,R(u,\xi)]+\frac{1}{2}[u, \prod [\xi, \eta]]-\frac{1}{2}[u,\prod [\xi, \eta]]+\\+2R(u, \prod [\xi, \eta])-\prod [u, [\xi,\eta]]+\prod [u, \prod [\xi,\eta]]   \Bigg\}=$\\

$=\prod [\xi, \eta]+\frac{1}{2}[\xi,[u,\eta]]+\frac{1}{2}\prod [\xi,\prod [u, \eta]]+2\prod [\xi, R(u, \eta)]-\frac{1}{2}\prod [\eta, [u, \xi]]-\\-\frac{1}{2}\prod [\eta, \prod [u, \xi]]-2\prod [\eta,R(u, \xi)]-\prod [u, [\xi, \eta]]=$\\

$=\prod [\xi, \eta]+\frac{1}{2}[\xi,[\eta, u]]-\frac{1}{2}\prod [\xi,\prod [ \eta, u]]+2\prod [\xi, R(u, \eta)]-\frac{1}{2}\prod [\eta, [\xi, u]]+\\+\frac{1}{2}\prod [\eta, \prod [\xi, u]]-2\prod [\eta,R(u, \xi)]$
\\
\\
where

$$\widetilde{a_{u}}(\xi, \eta)=-\frac{1}{2}\prod [\xi,\eta]-\frac{1}{4}\prod [[\xi, u],\eta]+\frac{1}{4}\prod [\prod [\xi, u],\eta]-\prod [R(u, \xi), \eta]+$$

$$+\frac{1}{4}\prod [[\eta, u],\xi]-\frac{1}{4}\prod [\prod [\eta, u], \xi]+\prod [R(u, \eta), \xi].$$\\

From this last equation it follows that:

$$\overset{1}{\nabla_{l}}a_{jk}^{i}\xi^{j}\eta^{k}\zeta^{l}=\frac{d}{dt}\widetilde{a_{\exp t\zeta}}(\xi, \eta)\arrowvert_{t=0}=-\frac{1}{4}\prod [[\xi, \zeta],\eta]+\frac{1}{4}\prod [\prod [\xi, \zeta],\eta]-\prod [R(\xi, \zeta), \eta]+$$

$$+\frac{1}{4}\prod [[\eta, \zeta],\xi]-\frac{1}{4}\prod [\prod [\eta, \zeta], \xi]+\prod [R(\eta, \zeta), \xi].$$\\

We obtain a result in conformity with proposition 2.1 and the relation $(3.2)$
in deed from the formulae $(3.2)$ it follows:

$\frac{1}{2}[b(\xi, \zeta, \eta)-b(\eta, \zeta, \xi)]=\frac{1}{2}\Bigg\{ -\frac{1}{2}\prod [[\xi, \zeta],\eta]+\frac{1}{2}\prod [\prod [\xi, \zeta],\eta]-2\prod [R(\xi, \zeta), \eta]+\\+\frac{1}{2}\prod [[\eta, \zeta], \xi]-\frac{1}{2}\prod [\prod [\eta, \zeta], \xi]+2\prod [R(\eta, \zeta),\xi]  \Bigg\}=$

$=-\frac{1}{4}\prod [[\xi, \zeta],\eta]+\frac{1}{4}\prod [\prod [\xi, \zeta], \eta]-\prod [R(\xi, \zeta),\eta]+\frac{1}{4}\prod [[\eta, \zeta], \xi]-$ 
$$-\frac{1}{4}\prod [\prod [\eta, \zeta], \xi]+\prod [R(\eta, \zeta),\xi].$$
\\

Therefore $$\overset{1}{\nabla_{l}}a_{jk}^{i}=b_{[j\arrowvert j \arrowvert k]}^{i}$$

\section{Computation of the tensor $d_{jklm}^{i}=\stackrel{2}{\nabla_{m}}b_{jkl}^{i}$}

 Denote $u\cdot R(\eta, \eta)\cdot u^{-1}$ by $R''_{u}(\eta, \eta)$. For the
computation of $d_{jklm}^{i}$ let us firstly compute  $R''_{u}(\eta, \eta)$.

The following proposition holds:\\

{\bf Proposition 5.1}\\

  $$R''_{u}(\eta, \eta)=R(\eta, \eta)+\prod [u, R(\eta, \eta)]+0(u, \eta^{2}) \eqno {(5.1)} .$$

The proof of this proposition, is from section 1 It is clear that
$\xi + \phi(\xi) \in Q$ and from section 4 $h''_{u}=h_{1}+[u, h_{1}]+0(u)$
where $h_{1} \in h$. Furthermore $\eta+R''_{u}(\eta, \eta) \in Q$ but
$ R''_{u}(\eta, \eta) \in h''_{u}$ that is why  $R''_{u}(\eta, \eta)$ can be 
represented as  $R''_{u}(\eta,\eta)=h_{1}+[u,h_{1}]+0(u)$, where 
$h_{1}=R''_{u}(\eta, \eta)-[u,R''_{u}(\eta, \eta)]+0(u)$. Let us write
 $\eta+R''_{u}(\eta, \eta) $ as :

$$ \eta+R''_{u}(\eta, \eta)=$$ $$=\Bigg\{(\eta+\prod [u,R''_{u}(\eta, \eta)])+( R''_{u}(\eta, \eta)-[u, R''_{u}(\eta, \eta) ])+([u, R''_{u}(\eta, \eta)]-\prod [u, R''_{u}(\eta, \eta)])  \Bigg\} $$ 
put $\eta+\prod [u, R''_{u}(\eta, \eta)]=\xi$ then\\
$\phi(\xi)= R''_{u}(\eta, \eta)-[u, R''_{u}(\eta, \eta)]+[u, R''_{u}(\eta, \eta)]-\prod [u, R''_{u}(\eta, \eta)]= R''_{u}(\eta, \eta)-\prod [u, R''_{u}(\eta, \eta)]+o(u) $

from the relation $(1.1)$ we have $\phi(\xi)=R(\xi, \xi)+S(\xi, \xi, \xi)+0(3)$

Therefore by comparing the term on the right hand sides of the last two relation,
we obtain:
$$R''_{u}(\eta, \eta)=R(\eta, \eta)+\prod [u, R(\eta, \eta)]+0(u, \eta^{2}).$$

Let  $\prod''_{u} :\mathfrak{G}  \longrightarrow V=T_{e}Q$ be the projection of
$\mathfrak{G}$ to $V$ parallel to $\mathfrak{h}''_{u}$. Then we obtain the equation
$$\xi+\widetilde{h}=\widetilde{\xi}+h_{1}+[u,h_{1}]$$

where $\xi,\widetilde{\xi} \in V $ and $\widetilde{h}, h_{1} \in \mathfrak{h}$
for the search of $\widetilde{\xi}=\widetilde{\xi}(\xi,u)$ we have

$$\xi+\widetilde{h}=\widetilde{\xi}+h_{1}+\prod [u,h_{1}]+([u,h_{1}]-\prod [u,h_{1}])$$

where $$\xi=\widetilde{\xi}+\prod [u,h_{1}]$$
$$\widetilde{h}=h_{1}+[u, h_{1}]-\prod [u, h_{1}]=h_{1}+\textrm{terms with u}. $$
From these two equalities we obtain

$$\widetilde{\xi}=\xi-\prod [u,\widetilde{h}]+0(u).$$
Hence
$$\prod''_{u}(\xi+\widetilde{h})=\xi-\prod [u,\widetilde{h} ]). \eqno {(5.2)} $$

We pass now to the computation of $d_{jklm}^{i}$.\\

From $(3.2)$ its follows that

$$b(\xi,\eta,\zeta)=-\frac{1}{2}\prod [[\xi,\eta],\zeta]+\frac{1}{2}\prod [\prod [\xi, \eta],\zeta]-2\prod [R(\xi,\eta),\zeta]$$
that is why
$$b''_{u}(\xi,\eta,\zeta)=-\frac{1}{2}\prod''_{u} [[\xi,\eta],\zeta]+\frac{1}{2}\prod''_{u} [\prod_{u} [\xi, \eta],\zeta]-2\prod''_{u} [R''_{u}(\xi,\eta),\zeta].$$

From $(5.2)$ it follows that\\

$$-\frac{1}{2}\prod''_{u}[[\xi,\eta],\zeta]=-\frac{1}{2}\prod [[\xi, \eta],\zeta]+\frac{1}{2}\prod [u,[[\xi, \eta], \zeta]]-\frac{1}{2}\prod [u,\prod [[\xi, \eta,], \zeta]]. \eqno {(5.3)} $$

Further more

$$\frac{1}{2}\prod''_{u} [\prod''_{u} [\xi, \eta],\zeta]=\frac{1}{2}\prod''_{u} [\prod [\xi, \eta], \zeta]-\frac{1}{2}\prod''_{u} [\prod [u,[\xi, \eta]], \zeta]+\frac{1}{2}\prod''_{u} [\prod [u, \prod [\xi , \eta]], \zeta]=$$

$$=\frac{1}{2}\prod [\prod [\xi, \eta], \zeta]-\frac{1}{2}\prod [u,[\prod[\xi, \eta], \zeta]]+$$

$$+\frac{1}{2}\prod [u, \prod [\prod [\xi , \eta], \zeta]]+o(u). \eqno {(5.4)}$$
\\
Finally from $(5.1)$ and $(5.2)$ it follows that:\\

$$-2\prod''_{u}[R''_{u}(\xi, \eta),\zeta]=-2\prod''_{u}[R_{u}(\xi, \eta),\zeta]-2\prod''_{u}[\prod [u, R(\xi, \eta),\zeta]$$
$$=-2\prod [R(\xi, \eta),\zeta]+2\prod [u, [R(\xi, \eta),\zeta ]]-2\prod [u,\prod [R(\xi, \eta), \zeta]]-$$
$$-2\prod [\prod [u, R(\xi,\eta),zeta]+o(4). \eqno {(5.5)}$$
\\
from (5.3), (5.4) and (5.5) it follows

$$d(\xi, \eta, \zeta, \tau)=\overset{2}{\nabla_{m}}b_{jkl}^{i}\arrowvert_{(e,e)}\xi^{j}\eta^{k}\zeta^{l}\tau^{m}=\frac{d}{dt}\Bigg(b''_{\exp t\tau}(\xi, \eta, \zeta)\Bigg)\arrowvert_{t=0}=$$

$$=\frac{1}{2}\prod [\tau, [[\xi, \eta], \zeta]]-\frac{1}{2}\prod [\tau, \prod [[\xi, \eta], \zeta]]-\frac{1}{2}\prod [\tau, [\prod [\xi, \eta ], \zeta ]]+\frac{1}{2}\prod [\tau, \prod [\prod [\xi, \eta], \zeta ]]-$$

$$-\frac{1}{2}\prod [\prod [\tau, [\xi, \eta ]],\zeta]+\frac{1}{2}\prod [\prod [\tau, \prod [\xi, \eta ]],\zeta]+2\prod [\tau,[R(\xi, \eta), \zeta]]-$$

$$-2\prod [\tau, \prod [R(\xi, \eta), \zeta]]-2\prod [\prod [\tau, R(\xi,\eta)],\zeta].  \eqno {(5.6)} $$
In the theory of 3-Webs \cite{ak, she1, she3} the following relation is known:
$$d_{jk[lm]}^{i}=-b_{jkp}^{i}a_{lm}^{p}.$$
Let us verify it:

$$\frac{1}{2}(d(\xi, \eta, \zeta, \tau)-d(\xi, \eta, \tau, \zeta))=\frac{1}{4}\prod [\tau, [[\xi, \eta],\zeta]]-\frac{1}{4}\prod [\zeta, [[\xi, \eta],\tau]]-\frac{1}{4}\prod [\tau,\prod [[\xi, \eta],\zeta]]+$$

$$+\frac{1}{4}\prod [\zeta, \prod [[\xi,\eta],\tau]]-\frac{1}{4}\prod [\tau, [\prod [\xi, \eta], \zeta]]+\frac{1}{4}\prod [\zeta, [\prod [\xi,\eta],\tau]]+\frac{1}{4}\prod [\tau, \prod [\prod [\xi, \eta],\zeta]]-$$

$$-\frac{1}{4}\prod [\zeta, \prod [\prod [\xi, \eta],\tau]]-\frac{1}{4}\prod [\prod [\tau, [\xi,\eta]],\zeta]+\frac{1}{4}\prod [\prod [\zeta,[\xi, \eta]], \tau]+$$

$$+\frac{1}{4}\prod [\prod [\tau, \prod [\xi, \tau]],\zeta]-\frac{1}{4}\prod [\prod [\zeta, \prod [\xi, \eta]],\tau]+\prod [\tau,[R(\xi,\eta),\zeta]]-\prod [\zeta, [R(\xi,\eta),\tau]]-$$

$$-\prod [\tau, \prod [R(\xi,\eta),\zeta]]+\prod [\zeta,\prod [R(\xi,\eta),\tau]]-\prod [\prod [\tau, R(\xi,\eta),]\zeta]+\prod [\prod [\zeta, R(\xi, \eta)],\tau]= $$\\
$$=-\frac{1}{4}\prod [[\xi,\eta],[\zeta,\tau]]+\frac{1}{4}\prod [\prod [\xi, \eta],[\zeta,\tau]]-\prod [R(\xi,\eta),[\zeta, \tau]].$$

In addition, considering that
$$[\zeta, \tau]=\prod [\zeta, \tau]+([\zeta, \tau]-\prod [\zeta,\tau]).$$

One obtain\\
$$\frac{1}{2}(d(\xi, \eta, \zeta, \tau)-d(\xi, \eta, \tau, \zeta))=-\frac{1}{4}\prod [[\xi, \eta],\prod [\zeta, \tau]]+\frac{1}{4}\prod [\prod [\xi, \eta], \prod [\zeta, \tau]]-$$
$$-\prod [R(\xi, \eta),\prod [\zeta, \tau]]$$

From relations $(3.1)$ and $(3.2)$ it follows that:\\

$$b(\xi, \eta, a(\zeta, \tau))=\frac{1}{2}b(\xi, \eta, \prod [\zeta, \tau])=-\frac{1}{4}\prod [[\xi, \eta],\prod [\zeta, \tau]]+\frac{1}{4}\prod [\prod [\xi, \eta], \prod [\zeta, \tau]]-$$
$$-\prod [R(\xi, \eta),\prod [\zeta, \tau]].$$
\\
Hence $d_{jk[lm]}^{i}=-b_{jkp}^{i}a_{lm}^{p}$

\section{Hexagonal loops}

The analytic hexagonal 3-Webs and their corresponding loops can be charaterise
by the following condition:
$$b_{(jkl)}^{i}=0$$

where $b(\xi, \eta, \zeta)=-\frac{1}{2}\prod[[\xi, \eta],\zeta]+\frac{1}{2}\prod[\prod [\xi, \eta],\zeta]-2\prod [R(\xi, \eta),\zeta]$

that is way, $b_{(jkl)}^{i}=0$ is equivalent to the following condition

$$\prod [R(\xi, \eta), \zeta]+\prod [R(\eta,\zeta),\xi]+\prod [R(\zeta, \xi),\eta]=0  \eqno {(6.1)}$$

which can be written as follows:

$$\underset{\xi \eta \zeta}{\sigma}\prod [R(\xi,\eta),\zeta]=0$$

where $\underset{\xi \eta \zeta}{\sigma}$ is the cyclic sum for $\xi, \eta, \zeta$

We have furthermore, for the hexagonal three-Webs the following relation 
$$d_{(jkl)m}^{i}=0.$$   Considering $(5.6)$ and $(6.1)$ one obtain

$$\underset{\xi \eta \zeta}{\sigma}\prod \Bigg\{[, \tau, [R(\xi,\eta),\zeta]]-[\prod [\tau, R(\xi, \eta)], \zeta]\Bigg\}=0.$$

where $\underset{\xi \eta \zeta}{\sigma}$ is the cyclic sum for $\xi, \eta, \zeta$

\end{document}